\documentstyle[12pt]{amsart}
\def\opn#1#2{\def#1{\operatorname{#2}}} 
\opn\chara{char}
\opn\length{\ell}
\opn\projdim{proj\,dim}
\opn\injdim{inj\,dim}
\opn\rank{rank}
\opn\depth{depth}
\opn\grade{grade}
\opn\height{height}
\opn\embdim{emb\,dim}
\opn\codepth{codepth}
\opn\codim{codim}

\opn\Tr{Tr}
\opn\bigrank{big\,rank}
\opn\superheight{superheight}\opn\lcm{lcm}
\opn\trdeg{tr\,deg}%
\opn\reg{reg}
\opn\ini{in}
\opn\div{div}
\opn\Div{Div}
\opn\cl{cl}
\opn\Cl{Cl}
\opn\Spec{Spec}
\opn\Supp{Supp}
\opn\supp{supp}
\opn\Sing{Sing}
\opn\Ass{Ass}
\opn\Ann{Ann}
\opn\Rad{Rad}
\opn\Soc{Soc}
\opn\Ker{Ker}
\opn\Coker{Coker}
\opn\Im{Im}
\opn\Hom{Hom}
\opn\Tor{Tor}
\opn\Ext{Ext}
\opn\End{End}
\opn\Aut{Aut}
\opn\id{id}

\opn\nat{nat}
\opn\pff{pf}
\opn\Pf{Pf}
\opn\GL{GL}
\opn\SL{SL}
\opn\mod{mod}
\opn\ord{ord}
\opn\aff{aff}
\opn\con{conv}
\opn\relint{relint}
\opn\st{st}
\opn\lk{lk}
\opn\cn{cn}
\opn\core{core}
\opn\vol{vol}
\opn\gr{gr}

\def\pot#1#2{#1[\kern-0.28ex[#2]\kern-0.28ex]}

\opn\dirlim{\underrightarrow{\lim}}
\opn\invlim{\underleftarrow{\lim}}
\def\Implies{\ifmmode\Longrightarrow \else
     \unskip${}\Longrightarrow{}$\ignorespaces\fi}
\def\implies{\ifmmode\Rightarrow \else
     \unskip${}\Rightarrow{}$\ignorespaces\fi}
\def\iff{\ifmmode\Longleftrightarrow \else
     \unskip${}\Longleftrightarrow{}$\ignorespaces\fi}

\let\:=\colon
\newtheorem{Theorem}{Theorem}[section]

\let\epsilon=\varepsilon
\let\phi=\varphi
\let\kappa=\varkappa
\def\qed{\ifhmode\textqed\fi
   \ifmmode\ifinner\quad\qedsymbol\else\dispqed\fi\fi}
\def\textqed{\unskip\nobreak\penalty50
    \hskip2em\hbox{}\nobreak\hfil\qedsymbol
    \parfillskip=0pt \finalhyphendemerits=0}
\def\dispqed{\rlap{\qquad\qedsymbol}}

\begin{document}
\title{Simultaneous resolution of singularities \\
(to appear in Proc. AMS)}
 \author{Steven Dale Cutkosky}
\thanks{partially supported by NSF. }
\date{}
 \maketitle
 \ \\ \\

\begin{abstract} We prove a local theorem on simultaneous resolution of singularities, which is valid
in all dimensions. This theorem  is proven in dimension 2 (and in all characteristics) by
Abhyankar in his book ``Ramification theoretic methods in algebraic geometry'' \cite{Ab2}. 
\end{abstract}

\section{Introduction}

Suppose that $f:V\rightarrow W$ is a generically finite morphism of normal varieties. 
If the characteristic of the ground field is zero, or if we are in char $p$ and dimension $\le 3$,
we know from resolution of singularities and resolution of indeterminancy of morphisms (\cite{H}, \cite{Ab3})  that
there are sequences of blowups of nonsingular subvarieties (monoidal transforms)
$V_1\rightarrow V$ and $W_1\rightarrow W$ such that $V_1$ and $W_1$ are nonsingular, and
$g:V_1\rightarrow W_1$ is a generically finite  morphism. 

A simultaneous resolution of $V$ and $W$ is a pair of nonsingular
$V_1$ and $W_1$ as above, such that $V_1\rightarrow W_1$ is a finite morphism.
Even if $V$ and $W$ have dimension two, a simultaneous resolution may not exist, even locally along a valuation,
 as shown by  an example of Abhyankar (Theorem 12 \cite{Ab4}).  

If we only require  that $W_1$ be nonsingular, but allow $V_1$ to be normal, we can take $V_1$ to be the normalization
of a resolution of singularities $W_1$ of $W$ in the function field of $V$, and obtain a  finite map
$V_1\rightarrow W_1$.  

However, the essential problem is to obtain a finite map such that $V_1$ is nonsingular.

In Theorem 1.1 we prove that given a fixed valuation $\nu$ of the characteristic zero function field $K$ of $V$, we can find a
nonsingular variety $V_1$, birationally dominating $V$, and a normal variety $W_1$ birationally dominating $W$, such that
$g:V_1\rightarrow W_1$ is finite in a neighborhood of the center of $\nu$. 
That is, if $p$ is the center of $\nu$ on $V_1$ and $q=g(p)$ is the center of $\nu$ on $W_1$,
then $\cal O_{W_1,q}$ lies below $\cal O_{V_1,p}$.

\begin{Theorem} 
Let $k$ be a field of characteristic zero, $L/k$  an  algebraic function field,
$K$ a finite algebraic extension of $L$, $\nu$ a  valuation of $K/k$, and $(R,M)$ a regular
 local ring with quotient field $K$, essentially of finite type over $k$, such that $\nu$ dominates $R$.
 Then for some sequence of monodial transforms $R\rightarrow R^*$ along $\nu$,
there exists a normal local  ring $S^*$ with quotient field $L$, essentially of finite type over $k$,
such that $R^*$ is the localization of the integral closure $T$ of $S^*$ in $K$ at a maximal ideal of $T$.
\end{Theorem}

We review notation from \cite{Ab2}.
A local domain $(R,M)$ is algebraic with ground field $k$ if $R$ is essentially of finite type over $k$.
Suppose that  $k, L, K$ are as in the statement of  Theorem 1.1. Suppose that $(S^*,N^)$ and $(R^*,M^*)$
 are normal algebraic local domains with ground field $k$ such that the quotient field of $S^*$ is $L$ and the
quotient field of  $R^*$ is $K$. We say that $S^*$ lies below $R^*$ if $R^*$ is a localization of the integral
closure $T$ of $S^*$ in $K$ at a maximal ideal of $T$. 
Suppose that $K$ is Galois over $L$ and $S^*$ lies below $R^*$. The splitting field $F^s(R^*/S^*)$
is the smallest field $K'$ between $L$ and $K$ such that $R^*$ is the only local ring in $K$ lying above $R^*\cap K'$. 
If $\nu$ is a valuation of $L$ with valuation ring $R_{\nu}$,  and  $\nu^*$ is an extension of $\nu$ to $K$
with valuation ring $R_{\nu^*}$, then $F^s(\nu^*/\nu)$ is defined to be $F^s(R_{\nu^*}/R_{\nu})$.

The following Theorem 1.2 follows immediately from our Theorem 1.1.

\begin{Theorem} Let $k$ be a field of characteristic zero, $L/k$  an $n$-dimensional algebraic function field,
$K$ a finite algebraic extension of $L$, $\nu$ a zero dimensional valuation of $K/k$, and $(R,M)$ a regular
algebraic local domain with quotient field $K$ and ground field $k$ such that $\nu$ has center $M$ in $R$.
 Then for some sequence of monoidal transforms $R\rightarrow R^*$ along $\nu$,
there exists an algebraic local domain $S^*$ with quotient field $L$ and ground field $k$ lying below $R^*$. 
\end{Theorem}

This Theorem is proved in dimension two (and in any characteristic) by  
Abhyankar  (Theorem 4.8 \cite{Ab2}). Our Theorem 2  is stated in the notation of Theorem 4.8  \cite{Ab2}.
 In the special case when the valuation $\nu$ has maximal rational rank n, 
Fu \cite{F} has shown that the conclusions of Theorem 2  hold. 

In Theorem 3 we  prove simultaneous resolution in a setting which is useful in the theory of resolution of singularities.

\begin{Theorem} Let $K$ be an $n$-dimensional algebraic function field with ground field $k$ of characteristic zero.
Let $\nu$ be a zero dimensional valuation of $K/k$. Suppose that $K^*$ is a Galois extension of $K$, and
$\nu^*$ is an extension of $\nu$ to $K^*$, $(R^*,M^*)$ and $(R,M)$ are $n$-dimensional normal algebraic local
domains with ground field $k$ and quotient field $K^*$ and $K$ respectively, such that $\nu^*$ has center $M^*$ in
$R^*$ and $R=R^*\cap K$. Suppose that $F^s(\nu^*/\nu)=F^s(R^*/R)$ and $R^s=R^*\cap F^s(\nu^*/\nu)$ is a
regular algebraic domain. Then there exists a sequence of monoidal transforms $R^s\rightarrow \overline R^s$ along $\nu^*$
such that there exists a regular algebraic domain $\overline R$ with ground field $k$ and quotient field $K$ lying below 
$\overline R^s$.
\end{Theorem}

In Theorem 4.9 \cite{Ab2}, Abhyankar
proves local uniformization along a valuation in dimension 2, when the  ground field $k$ is 
algebraically closed  of characteristic zero.
 The proof in fact shows more. It proves that 
(over an algebraically closed ground field of characteristic zero), ``embedded local uniformization in algebraic local
domains of dimension $n$'' + ``the conclusions of our Theorem 2 (or Theorem 3) in dimension $n$''
 implies ``local uniformization in dimension $n$''.

Our theorems are  an application of our theorem on 
monomialization of generically finite  morphisms along a valuation (Theorem A \cite{C2}) and the general theory
developed in \cite{Ab2}. We thank Professor
Abhyankar for pointing out to us that simultaneous resolution is related to our monomialization theorem.
Here we state the monomialization theorem from our paper \cite{C2}, which we use in this paper.

\begin{Theorem} (Theorem A \cite{C2} - Monomialization)
 Suppose that $R \subset S$ are  excellent regular local rings such that $\text{dim}(R)=\text{dim}(S)$,
 containing a field $k$ of
characteristic zero, such that the quotient field $K$ of $S$ is a finite extension of the quotient field $J$
of $R$.

Let $V$ be a valuation ring of $K$ which dominates $S$. 
Suppose that if $m_V$ is the maximal
ideal of $V$, and $p^*=m_V\cap S$, then $(S/p^*)_{p^*}$ is a finitely generated field extension of $k$.
 Then there
exist sequences of  monoidal transforms (blow ups of regular primes)
 $R \rightarrow R'$ and $S \rightarrow S'$
such that $V$ dominates $S'$, $S'$ dominates $R'$ and 
there are regular parameters $(x_1, .... ,x_n)$
in $R'$,  $(y_1, ... ,y_n)$ in $S'$, units $\delta_1,\ldots,\delta_n\in S'$ and a matrix $(a_{ij})$ of
nonnegative integers such that  $\text{Det}(a_{ij}) \ne 0$ and
$$
\begin{array}{ll}
x_1 &= y_1^{a_{11}} ..... y_n^{a_{1n}}\delta_1\\
\vdots&\\   
x_n &= y_1^{a_{n1}} ..... y_n^{a_{nn}}\delta_n.
\end{array}
$$

\end{Theorem}

Theorem A is used in \cite{C2} to prove a factorization theorem for birational morphisms. This theorem is proved in 
dimension 3 in our paper \cite{C1}.

\section{Simultaneous resolution}

\subsection{Resolution of Singularities}

Hironaka's theorems on resolution, the ``fundamental theorems'' and ``main theorems'' stated in chapter I of \cite{H},
apply to ``algebraic schemes''. An algebraic scheme is defined (on page 162 of \cite{H}) to be a separated scheme of finite type
over $\text{spec}(S)$, with $S$ a local ring in $\beta$. The class $\beta$ of local rings is defined on page 161 of \cite{H}.
$S\in \beta$ if

\begin{enumerate}
\item $S$ is a (Noetherian) local ring with residue field of characteristic 0.
\item If $\hat S$ denotes the completion of the local ring $S$, then for every $S$-algebra of finite type $A$, the singular
locus of $\text{spec}(A\otimes_S\hat S)$ is the preimage of that of $\text{spec}(A)$ under the canonical map of the
first spectrum into the second.
\end{enumerate}

\begin{Theorem} Suppose that $S$ is an excellent local ring containing a field of characteristic 0. Then $S\in \beta$.
Thus $\text{spec}(S)$ is an ``algebraic scheme''.
\end{Theorem}

\begin{pf}
For a Noetherian ring $B$, let $\text{reg}(B)$ denote the subset of $\text{spec}(B)$ of primes $p$ such that $B_p$ is
a regular local ring,  $\text{sing}(B)$  denote the subset of $\text{spec}(B)$ of primes $p$ such that $B_p$ is
not a regular local ring.

Suppose that $A$ is an $S$-algebra of finite type. $S\rightarrow \hat S$ 
is faithfully flat and regular, since $S$ is a local $G$-ring. 
The natural map $\phi:A\rightarrow A\otimes_S\hat S$ is then faithfully flat and regular
by Lemma 4 (33.E) \cite{M}. 
$(\phi^*)^{-1}(\text{reg}(A)) = \text{reg}(A\otimes_S\hat S)$
by Theorem 51 (21.D) \cite{M}. Thus 
$(\phi^*)^{-1}(\text{sing}(A)) = \text{sing}(A\otimes_S\hat S)$.
\end{pf}

Suppose that $(R,m)$ is a  local domain, with maximal ideal $m$, and that $P\subset  R$
is a prime ideal, such that $R/P$ is regular. Suppose that $0\ne f\in P$, and $m_1$ is a prime ideal in $R[\frac{P}{f}]$
such that $m_1\cap R = m$. Set 
$R_1 =(R[\frac{P}{f}])_{m_1}$.
$R_1$ (or $R\rightarrow R_1$) is called a monoidal transform of $R$. If $P=m$, then $R_1$ is called a
quadratic transform.

We will say that a valuation $\nu$ of the quotient field of $R$ dominates $R$ if $\nu$ has nonnegative value on $R$, and 
$m = \{f\in R | \nu(f) >0\}$. A monoidal transform $R\rightarrow R_1$ is along $\nu$ if $\nu$
dominates $R$ and $R_1$.

\begin{Theorem} Suppose that $R$, $S$ are excellent local domains containing a field $k$ of characteristic zero
such that $S$ dominates $R$ and $S$ is regular. Let $\nu$ be a valuation of the quotient field $K$ of $S$ that dominates $S$,
$R\rightarrow R_1$ a monoidal transform such that $\nu$ dominates $R_1$. 
$R_1$ is a local ring on $X=\text{Proj}(\bigoplus_{n\ge0}p^n)$ for some prime $p\subset R$.
Let
$$
U=\{Q\in\text{spec}(S):pS_Q\text{ is invertible }\}
$$
an open subset of $\text{spec}(S)$. Then there exists a projective morphism  $f:Y\rightarrow \text{spec}(S)$
which is a product of monoidal transforms
such that if $S_1$ is the local ring of $Y$ dominated by $\nu$, then   $S_1$ dominates $R_1$, and
$(f)^{-1}(U)\rightarrow U$ is an isomorphism.
\end{Theorem}
\begin{pf}
Since $S$ is a UFD, we can write $pS=gI$, where $g\in S$, $I\subset S$ has height $\ge 2$. Then $U=\text{spec}(S)-V(I)$.
By Main Theorem II(N) \cite{H}, there exists a sequence of monoidal transforms $\pi:Y\rightarrow\text{spec}(S)$
such that $I\cal O_Y$ is invertible, and $\pi^{-1}(U)\rightarrow U$ is an isomorphism. 
Let $S_1$ be the local ring of the center of $\nu$ on $Y$.
We have $pS_1 = hS_1$ for some $h\in p$. Hence $R[\frac{p}{h}]\subset S_1$, and since $\nu$ dominates $S_1$,
$R_1$ is the localization of $R[\frac{p}{h}]$ which is dominated by $S_1$.
\end{pf}

\begin{Theorem} Suppose that $R$ is an excellent regular local domain containing a field of characteristic zero,
with quotient field $K$. Let $\nu$ be a valuation of $K$ dominating $R$. Suppose that $f\in K$ is such that $\nu(f)\ge0$.
Then there exists a sequence of monoidal transforms along $\nu$
$$
R\rightarrow R_1\rightarrow\cdots\rightarrow R_n
$$
such that $f\in R_n$.
\end{Theorem}

\begin{pf} Write $f=\frac{a}{b}$ with $a,b\in R$. By  Main Theorem II(N) \cite{H} applied to the ideal $I=(a,b)$ in $R$,
there exists a sequence of monoidal transforms along $\nu$,
 $R\rightarrow R_n$ such that $IR_n=\alpha R_n$ is a principal ideal. There exist
$c,d,u_1,u_2$ in $R_n$ such that $a=c\alpha,b=d\alpha,\alpha=u_1a+u_2b$. Then
$u_1c+u_2d=1$, so that $cR_n+dR_n=R_n$, and one of $c$ or $d$ is a unit in $R_n$. If $c$ is a unit, then
$0\le\nu(f)=\nu(\frac{c}{d})=\nu(c)-\nu(d)$ implies $\nu(d)=0$, and since $\nu$ dominates $R_n$,  $d$ is a unit and $f\in R_n$.
\end{pf}

Suppose that $Y$ is an algebraic scheme, $X$, $D$ are subschemes of $Y$. Suppose that $g:Y'\rightarrow Y$,
$f:X'\rightarrow X$ are the monoidal transforms of $Y$ and $X$ with center $D$ and $D\cap X$ respectively. Then there
exists a unique isomorphism of $X'$ to a subscheme $X''$ of $Y'$ such that $g$ induces $f$ (cf. chapter 0, section 2 \cite{H}).
$X''$ is called the strict transform of $X$ by the monoidal transform $g$.

\begin{Theorem} 
Let $R$ be an excellent regular local ring, containing a field of characteristic zero. Let $W\subset \text{spec}(R)$ be an integral
subscheme, $V\subset \text{spec}(R)$ be the singular locus of $W$. Then there exists a sequence of monoidal transforms
$f:X\rightarrow\text{spec}(R)$ such that the strict transform of $W$ is nonsingular in $X$, and $f$ is an isomorphism over 
$\text{spec}(R)-V$.
\end{Theorem}
\begin{pf} This is immediate from Theorem $I_2^{N,n}$ \cite{H}.
\end{pf}

\subsection{Proof of Theorem 1.1}
Let $V$ be the valuation ring of $\nu$ with maximal ideal $m_{\nu}$.

We will first prove the theorem with the assumption that $\text{trdeg}_k V/m_{\nu} = 0$.
Let $q_1,\ldots,q_n$ be a transcendence basis of $L/k$. After replacing $q_i$ by $\frac{1}{q_i}$ if necessary, we
may assume that $\nu(q_i)\ge 0$ for all $i$. By Theorem 2.3, after possibly replacing $R$ by a sequence of monoidal
tranforms, we may assume that $q_1,\ldots, q_n\in R$. 

Let $T$ be the integral closure of $k[q_1,\ldots,q_n]$ in $L$. $T\subset R$. Set $N'=M\cap T$, $S=T_{N'}$,
$N=N'T_{N'}$. $\text{trdeg}_k S/N\le \text{trdeg}_k R/M \le \text{trdeg}_k V/m_{\nu} =0$. Thus
$\text{dim }S = \text{trdeg}_k L - \text{trdeg}_k S/N = n$.

By Theorems 2.4 and 2.2, we can perform a sequence of monoidal transforms  $S\rightarrow S_1$ and 
$R\rightarrow R_1$  so that $\nu$ dominates $R_1$, $R_1$ dominates $S_1$,
$R_1$ and $S_1$ are regular, and $\text{dim }S_1 = \text{dim }R_1 = n$.

By Theorem 1.4, we can perform sequences of monoidal transforms  $R_1\rightarrow R^*$ and $S_1\rightarrow S^*$
so that $\nu$ dominates $R^*$, $R^*$ dominates $S^*$, $R^*$ 
has regular parameters $(y_1,\ldots,y_n)$ and $S^*$ has regular parameters $(x_1,\ldots,x_n)$, such that there are 
 units $\delta_1,\ldots,\delta_n\in R^*$ and a matrix $(a_{ij})$ of
nonnegative integers such that  $\text{Det}(a_{ij}) \ne 0$ and
$$
\begin{array}{ll}
x_1 &= y_1^{a_{11}} ..... y_n^{a_{1n}}\delta_1\\
\vdots&\\   
x_n &= y_1^{a_{n1}} ..... y_n^{a_{nn}}\delta_n.
\end{array}
$$ 
By permuting variables, we may asume that $D=\text{Det}(a_{ij}) > 0$. Let $(b_{ij})$ be the adjoint matrix of 
$(a_{ij})$. We have that 
$$
x_1^{b_{j1}}\cdots x_n^{b_{jn}}=y_j^D\epsilon_j
$$
for $1\le j\le n$, where $\epsilon_j$ are units in $R^*$. Since these elements are in $R^*\cap L$, we have the existence
of a normal algebraic local domain $S^*$ with ground field $k$ and quotient field $L$ lying below $R^*$ by
Proposition 3.12 \cite{Ab2}.

Now we will prove the theorem for $\text{trdeg}_k V/m_{\nu}$ arbitrary. 
Let $\nu'$ be the restriction of $\nu$ to $L$. The valuation ring of $\nu$ is $V' = V\cap L$, with maximal ideal 
$m_{\nu'} = m_{\nu}\cap L$. 
$$
\mbox{trdeg}_k V/m_{\nu} = \mbox{trdeg}_k V'/m_{\nu'} \le \mbox{trdeg}_k L <\infty
$$
by Proposition 2.46 \cite{Ab2}. We can lift a transcendence basis of $V'/m_{\nu'}$ over $k$ to $f_1,\ldots f_m\in L$.
$f_1,\ldots, f_m$ are algebraically independent over $k$. Set $k' = k(f_1,\ldots,f_m)$. By Theorem 2.3 we can
perform a sequence of monoidal transforms $R\rightarrow R'$ along $\nu$ so that $f_1,\ldots, f_m\in R'$.
$k'\subset R'$ since $f_1,\ldots, f_m$ are algebraically independent over $k$ and $\nu$ dominates $R'$. 

We have that $L$ is a field of algebraic functions over $k'$, $R'$ is essentially of finite type over $k'$
and $V$ is a valuation of $K/k'$ such that $\mbox{trdeg}_{k'} V/m_{\nu} = 0$. By the first part of the proof
there exists a 
 sequence of monodial transforms $R'\rightarrow R^*$ along $\nu$,
 a normal local  ring $S^*$ with quotient field $L$, essentially of finite type over $k'$,
such that $R^*$ is the localization of the integral closure $T$ of $S^*$ in $K$ at a maximal ideal of $T$.
$S^*$ is essentially of finite type over $k$ since $k'$ is essentially of finite type over $k$.

\subsection{Proof of Theorem 1.3}

By Theorem 1.2, there exists a monoidal transform sequence $R^s\rightarrow \overline R^s$ along $\nu^*$ such that 
there exists a normal algebraic domain $\overline R$ in $K$ lying below $\overline R^s$. Let $\overline M$,
 $\overline M^s$ be the respective maximal ideals in $\overline R$ and and $\overline R^s$. Let 
$(\overline R^*,\overline M^*)$ be the local ring in $K^*$ lying above $\overline R^s$ such that $\nu^*$
has center $\overline M^*$ in $\overline R^*$.
$$
F^s(\nu^*/\nu)=F^s(R^*/R)\subset F^s(\overline R^*/\overline R)\subset F^s(\nu^*/\nu)
$$
implies $F^s(\overline R^*/\overline R)=F^s(\nu^*/\nu)$.
Thus $\overline R^s/\overline M^s=\overline R/\overline M$, $\overline M \,\overline R^s = \overline M^s$
 and $\overline R$ is regular by
Theorem 1.47 \cite{Ab2} and Proposition 3.18B \cite{Ab2}.

\ \\ \\
\noindent
S. Dale Cutkosky, Department of Mathematics, University of
Missouri\\
Columbia, MO 65211\\
dale@@cutkosky.math.missouri.edu

\begin{thebibliography}{99}

\bibitem{Ab1} Abhyankar, S., Local uniformization on algebraic surfaces over ground fields of
characteristic $p\ne 0$, Annals of Math, 63 (1956), 491-526.
\bibitem{Ab2} Abhyankar, S., Ramification theoretic methods in algebraic geometry, Princeton 
University Press, 1959.
\bibitem{Ab3} Abhyankar, S., Resolution of singularities of embedded algebraic surfaces, Academic Press, New York, 1966.
\bibitem{Ab4} Abhyankar, S., Simultaneous resolution for algebraic surfaces, Amer. J. Math 78 (1956), 761-790.
 \bibitem{C1} Cutkosky, S.D.,
Local factorization of birational maps, Advances in Math. 132, (1997), 167-315.
\bibitem{C2} Cutkosky, S.D., Local factorization and monomialization of morphisms, to appear in Asterisque.
\bibitem{F} Fu, D., Local weak simultaneous resolution for high rational ranks, Purdue Univ. Thesis 1996.
\bibitem{H} Hironaka, H., Resolution of singularities of an algebraic variety over a field of
characteristic zero, Annals of Math, 79 (1964), 109-326.
\bibitem{M} Matsumura, H, Commutative Algebra, second edition, Benjamin/Cummings, Reading Mass, 1980.
\end{thebibliography}
\end{document}